\documentclass[11pt,a4paper]{article}
\usepackage{amsmath}
\usepackage{amsfonts}
\usepackage{amssymb}
\usepackage{amsopn,amsthm}
\usepackage[a4paper]{geometry}
 \usepackage{hyperref}
\usepackage[english]{babel}
\usepackage{cite}
\usepackage{fancyhdr}

\newtheorem{theorem}{Theorem}[section]
\newtheorem{lemma}[theorem]{Lemma}
\newtheorem{proposition}[theorem]{Proposition}

\theoremstyle{definition}
\newtheorem{definition}[theorem]{Definition}
\newtheorem{example}[theorem]{Example}

\newtheorem{remark}[theorem]{Remark}

\newcommand{\W}{\mathcal{W}}\newcommand{\Ric}{\mathrm{Ric}}

\newcommand{\grad}{{\mathrm{grad}}}

 \renewcommand{\Im}{{\operatorname {Im}}}
\renewcommand{\Re}{{\operatorname {Re}}}
 \renewcommand{\rho}{ {\varrho}}

 \renewcommand{\epsilon}{ {\varepsilon}}

\newcommand{\abs}[1]{|#1|}

\def\T{{\mathcal T}}
\def \e {{\varepsilon}}
\def \g {{\gamma}}
\def \ga {{\gamma}}

\def \a {{\alpha}}
\def \b {{\beta}}
\def \d {{\delta}}
\def \s {{\sigma}}
\def \la {{\lambda}}
\def \O {{\Omega}}
\def \m {{\mu}}
\def \L {{\Lambda}}

\def \z {{\zeta}}

\def \r {{\varrho}}
\def \endp { $\hfill \square$ \vskip 8 pt}

\def \R {{\mathbb {R}}}
\def \H {{\mathcal {H}}}
\def \C {{\mathbb {C}}}

\def \N {{\mathbb N}}
\def \t {{\tau}}
\def \wt {\widetilde}
\def \ol {\overline}

\def \p {{\partial}}

\renewcommand{\T}{{T^{1,0} }}
\newcommand{\WW}{{\overline W }}
\newcommand{\VV}{{\overline V }}
\newcommand{\ZZ}{{\overline Z }}
\newcommand{\E}{{\mathcal{E} }}
\newcommand{\w}{{\overline w} }
\newcommand{\ww}{{\overline w }}

\newcommand{\1} {{\bar{1} }}
\newcommand{\2}{{\bar{2} }}
\newcommand{\EE}{{\overline{E} }}
\newcommand{\UU}{{\overline{U} }}
\newcommand{\TT}{{T^{0,1} }}
\newcommand{\zz}{{\overline{z}}}
\renewcommand{\ss}{{\overline\sigma }}
\renewcommand{\aa}{{\overline{\alpha}}}
\newcommand{\bb}{{\overline{\beta}}}
\renewcommand{\l}{{{\lambda}}}
\renewcommand{\ll}{{\overline{\lambda}}}
\newcommand{\kk}{{\overline{k}}}
\newcommand{\rr}{{\overline{\rho}}}

\renewcommand{\bar}{\overline}
\newcommand{\D}{{\nabla}}
\renewcommand{\gg}{{\overline{\gamma}}}
\newcommand{\mm}{{\overline{\mu}}}
\newcommand{\barint}
{\rule[.036in]{.12in}{.009in}\kern-.16in \displaystyle\int}

\renewcommand{\theta}{\vartheta}

\begin{document}
\title{Pseudohermitian invariants and classification of CR~mappings in
generalized ellipsoids}

 \author{Roberto Monti and Daniele Morbidelli}


\date{}
 \maketitle


\begin{abstract}
Given  a strictly pseudoconvex  hypersurface $M\subset \C^{n+1}$,
we discuss the problem of classifying all local CR diffeomorphisms between open subsets  $N, N'
\subset M$.  Our method exploits the  Tanaka--Webster  pseudohermitian
invariants of a   contact form $\theta$ on $M$,  their
transformation formulae, and the
Chern--Moser invariants.
Our
main application concerns a class of generalized ellipsoids where we
classify all local CR mappings.
\end{abstract}

\setcounter{section}{-1}
\section{Introduction}

In this paper, we give a contribution to the problem of classifying     local CR
mappings between
real hypersurfaces in
$\C^{n+1}$.  Namely, given a surface $M:=b\Omega$, where
$\Omega\subset\C^{n+1}$ is a smooth open set, we consider the problem of
classifying all CR mappings $f: N\to N'$, where $N$ and $ N'$ are open subsets
of $M$.
The question  is rather natural, because
biholomorphic mappings of $\Omega$ that extend smoothly to the boundary
define CR mappings on $M$.

Our approach is mainly based on CR differential geometry of   strictly
pseudoconvex ma\-ni\-folds.
We fix a contact form $\theta$ on $M$,
we calculate the Tanaka--Webster invariants (see \cite{Tanaka62,W77}) and we
exploit Lee's transformation formulae,  see \cite{L}.
The idea is reminiscent of  known
techniques in the  study of  conformal mappings in  Riemannian manifolds, see
\cite{LF,Kuhnel,Iwaniec}.
 Our  point of view is
  described in Section \ref{ladue}, in the
setting of CR surfaces in $\C^{n+1}$,
$n\ge 2$.


We also exploit the connection between
pseudohermitian invariants and the classical
Cartan--Chern--Moser CR invariants, see  \cite{CM}.
In particular, we
introduce a new  Chern-invariant cone bundle.
Namely,
starting from the Chern tensor, we define   a subset $\H :  ={\bigcup_{P\in
M}\H_P}$ of the holomorphic tangent bundle which is
preserved by CR mappings.
The definition of $\H$ is given in Section \ref{ladue}.
We believe that the study  of this cone bundle    may be of some  interest in similar or related situations.

These ideas are applied to  the model given by a \emph{generalized ellipsoid}
\begin{equation}\label{cinquesei}
M:=b E,\quad E=\big\{ z\in \C^{n+1}: \abs{z_1}^{2m_1}+\cdots
+\abs{z_{s-1}}^{2m_{s-1}} +\abs{z_s}^2=1\big\},
\end{equation}
 where $z_1,z_2,\dots ,z_s $ are groups of variables and the numbers $m_j$
satisfy suitable hypoteses.
 The automorphism group of  $E$ is  studied in \cite{Sunada} and the
model is considered also in   \cite{Kodama}.  In \cite{DS} the authors
prove that all local CR mappings of $b E$ extend to global biholomorphic
mappings of
$E$,
under suitable hypotheses on the dimension of the groups of variables $z_j$.

In the present paper, we study local  CR
mappings on $b E$.
We recover both the
results in \cite{Sunada} and \cite{DS} on the model \eqref{cinquesei}.
Our arguments are completely different and new.
The statement
of our classification result for CR mappings on generalized ellipsoids is
contained in   Section
\ref{scheletro}, Theorem \ref{principale}.  Subsections \ref{cierre} and
\ref{imparo}
 contain the  computation of the pseudohermitian and Chern--Moser invariants in
our model.  Section \ref{SFRT} is devoted to the computation of the CR factor
of a  mapping and   to the classification of  CR mappings which are ``Levi-isometric''.

We mainly use differential geometric
arguments, which  require a certain computational effort. On the other hand,
they provide a good understanding
of the geometry of the manifold $M$. Several
other  strategies are available in the study of CR
mappings.
  For a complete account, we refer the reader to the monograph
\cite{BERo}, where  \emph{Segre varieties},  \emph{infinitesimal
CR mappings} and other tools are widely discussed.

\section{Chern-invariant cones}\label{ladue}
\setcounter{equation}{0}

Let $M\subset\C^{n+1}$ be a
strictly pseudoconvex real hypersurface.  Fix on $M$ a
contact form
$\theta$ and let $L: =-id\theta$ denote the Levi form. For a fixed
frame  of holomorphic vector fields $Z_\a$, $\a=1,\dots, n$, let
 $h_{\a\bb}=
L (Z_\a, Z_\bb)$.

A diffeomorphism  $f:N\to N'$
between open subsets $N $  and $ N'$ of $M$
is by definition a \emph{CR mapping}
if  $f^*\theta =
\lambda \theta$ for some
function $\lambda>0$ on $N$
and the tangent mapping $f_*$ preserves the complex structure. We call $\lambda$
the \emph{CR factor} of $f$.
Observe  that CR mappings preserve orthogonality with respect to
the Levi form:
\begin{equation}
 \label{caubo}
L(f_*Z, f_*\ol{W})= \lambda L(Z, \ol{W}) \qquad \text{for all $ Z, W\in \T N $,}
\end{equation}
where $\T N\subset\C TN$ denotes the holomorphic tangent bundle.
Observe also that
$L(f_*Z, f_*W)=L(Z, W)=0$ for all $Z,W\in \T N$.
Therefore CR mappings are somewhat similar to conformal mappings in Riemannian
manifolds
and \eqref{caubo} is an overdetermined system analogous to the system satisfied
by conformal mappings in the Riemannian setting.

Let us denote by $u =\lambda^{-1}$ the inverse of the CR factor.
Given a contact form $\wt\theta= u^{-1}\theta$,
the pseudohermitian Ricci curvature $R_{\a\bb}$ and the preusohermitian torsion
$A_{\a\b}$ transform according to Lee's formulae \cite{L}:
\begin{equation}
 \label{cari}
  \wt R_{\a \bar \b} = R_{\a\bar\b} +\frac{n+2}{2u}\Big\{
  u_{,\a\bar\b}+u_{,\bar\b\a} -\frac{2}{u} u_{,\a} u_{,\bar\b}\Big\}
  +\frac{1}{2u}
  \Big\{
     \Delta u - \frac{2(n+2)}{u} |\D u|^2\Big\} h_{\a\bar\b},
\end{equation}
and \begin{equation}
 \label{tors}
  \wt A_{\a\b} = A_{\a\b} - \frac{i}{u} u_{,\a\b}.
\end{equation}
We refer the reader to
\cite{Tanaka75,W77,DT} for definition and basic properties of these
tensors.
Here, $u_{,\a\bb} = \D_\bb\D_\a u$ denote second order covariant derivatives
with   respect to the Webster connection
$\D$, while $\Delta u = u_{,\g}{}^\g + u_{,\gg}{}^\gg$ and $\abs{\nabla u}^2 = u_\g
u^\g$.  As usual, we raise an  lower indices   by $h^{\a\bb}$, the
matrix defined by $ h^{\a\bb}h_{\g\bb}=\d_\g^\a$, so that
$u_,{}^\g= h^{\g\bb}u_{,\bb}$ and $ u_{,\g}{}^\g = h^{\g\bb} u_{,\g\bb} $.
Here and henceforth, we omit   summation on repeated indices.
For future reference, recall the contracted  version of  \eqref{cari}
\begin{equation}
 \label{change4}
\begin{split}
 \frac {1}{u} \wt R =  R+ \frac{n+1}{u}
  \left\{ \Delta u
 -\frac{n+2}{u} |\D u|^2 \right\},
\end{split}
\end{equation}
where $R:= h^{\a\bb}R_{\a\bb}= R_\a{}^\a$ is the pseudohermitian scalar
curvature.

Formulae \eqref{cari} and \eqref{tors} are relevant in the study of
the CR Yamabe problem, see \cite{JL}. The Riemannian version of
\eqref{cari} is also important in some regularity questions for conformal
mappings,
see \cite{LF}.

Equations \eqref{cari} and \eqref{tors} form a system  of nonlinear PDEs for
the inverse of the CR factor $u$. This system also involves $f$, in a
way that becomes clear in the coordinate-free notation:
\begin{equation}\label{ochetta}
\begin{aligned}
   \Ric (f_*Z, f_*\WW) &= \Ric(Z, \WW) +\frac{n+2}{2u}\Big\{
  \D^2 u (Z, \WW)+\D^2 u (\WW, Z)  -\frac{2}{u} Zu\WW u \Big\}
\\&\qquad  +\frac{1}{2u}
  \Big\{
     \Delta u - \frac{2(n+2)}{u} |\D u|^2\Big\} L(Z, \WW),
\end{aligned}\end{equation}
and
\begin{equation}\label{pablo}
 A(f_*Z,f_* W) = A(Z, W) - \frac{i}{u} \D^2 u(Z, W),
\end{equation}
 for
any $Z, W\in \T N$.
The function   $f$ appears  through the
geometric terms with $\Ric$ and $A$ in the left-hand side of \eqref{ochetta}
and \eqref{pablo}.

When $M$ has dimension $2n+1\ge 5$, i.e.~$n\geq 2$,
the Chern tensor $S_{\a\bb\l\mm}$ introduced in \cite{CM}  is a nontrivial
relative CR
invariant which satisfies
\begin{equation}
 \label{chevi}
  \wt S_{\a\bb\la\bar\mu}  = \frac{1}{u} S_{\a\bb\la\bar\mu},
\end{equation}
see  \cite{W00}.
The tensor $S_{\a\bb\l\mm}$ can be expressed in terms of the pseudohermitian
curvatures by means of the \emph{Webster's formula},
 see \cite{W77,DT},
\begin{equation}\begin{split}\label{web}
   S_{\a\bb\la\bar\mu}  = R_{\a\bb\la\bar\mu } & -
         \frac{1}{n+2} \big\{h_{\a\bb}R_{\la\bar\mu} +
       h_{\la\bb}R_{\a\bar\mu}+h_{\a\bar\mu}R_{\la\bb}
         +h_{\la\bar\mu}R_{\a\bb} \big\}
   \\&
+\frac{R}{(n+1)(n+2)}\big\{h_{\a\bb}h_{\la\bar\mu}+
    h_{\la\bb}h_{\a\bar\mu} \big\}.
\end{split}
\end{equation}

Let us introduce a  Chern-invariant cone bundle $\H\subset \T M$ which is
preserved by CR mappings.
Namely, let   $\H:=\displaystyle{\bigcup_{P\in { M}}}\H_P $,
where
\begin{equation}
\begin{split}
  \label{accadi}
  \H_P : = \Big \{ U\in T_P^{1,0} M  : R(U, \VV, & Z, \WW) = 0
   \text{ for all } V, Z, W\in  T_P^{1,0} M \text{ such that}
    \\&
    L(U, \VV)= L(U, \WW)= L(Z, \VV)= L(Z, \WW)=0 \Big \}.
\end{split}
\end{equation}
The set $\H_P$ is a cone in the vector space  $T_P^{1,0} M$.  If
$U,V,Z,W\in  T_P^{1,0} M $ satisfy the orthogonality relations in
\eqref{accadi}, then we have $ R(U, \VV, Z, \WW)=S(U, \VV, Z, \WW)$, by \eqref{web}.
 In particular, $\H$ does not depend on $\theta$ and, moreover,  any
CR mapping $f:N\to N'$
satisfies
\begin{equation}\label{accado}
 f_*(\H_P)= \H_{f(P)}  \qquad\text{ for any } P \in N.
\end{equation}
In general, the set $\H_P$ is not closed under addition. If $\theta$
is any  contact form on the standard sphere or on the Siegel domain, then $\H= \T
M$ and \eqref{accado} carries no information. On the other hand,
if the Chern tensor has a nontrivial structure, then $\H$ can
provide some useful  information on how  a
CR mapping $f$ transforms the holomorphic tangent space.

In the case of the generalized  ellipsoids \eqref{cinquesei},
 $\H$ can be decomposed as $\H=\mathcal{G} \oplus \mathcal{G}^\perp$,
where $ \mathcal{G}$ and $ \mathcal{G}^\perp$ are orthogonal
with respect to the Levi form.
This  is discussed in Section \ref{imparo}.
A cone similar to $\H$ was used in the Riemannian setting in
\cite{Morbidelli}. In that case, the cone bundle   is  related to
umbilical surfaces in the underlying structure. It could be of  interest to
understand whether also  in the CR setting the cone $\H$ is related to
significant geometric objects.

\section{The generalized ellipsoid model: main result and skeleton of the proof}
\setcounter{equation}{0}

In this section,  we state the classification theorem for
CR mappings on generalized ellipsoids and we indicate the scheme of the proof.
\label{scheletro}
Let
\[
p(z,\zz): = |z_1|^{2m_1}+\cdots +\abs{z_{s-1}}^{2m_{s-1}}+ \abs{z_s}^2,
\]
where $z=(z_1, \dots,  z_s)\in \C^{n_1}\times \cdots\times\C^{n_s}=\C^n$.
We assume that the integers $m_j,n_j$, $j=1,...,s$, satisfy
\begin{equation}
   \label{ipo}
\begin{cases}
  m_j > 1  \textrm{ and } n_j\ge 2, & \textrm{ if } 1\le j\le s-1,
\\
  n_s\ge 0. &\end{cases}
\end{equation}
We denote  by $z^\a$ the $\a$th variable in $\C^n$
and   we partition the indexes $\{1,...,n\}$ into the following
sets
\[
    I_1 = \{1,...,n_1\},\quad
    I_2 = \{n_1+1,...,n_1+n_2\},\quad
    ...\quad
    I_s = \{n_1+\cdots +n_{s-1}+1,...,n\},
\]
so that $\abs{z_j}^2=\sum_{\a\in I_j}\abs{z^\a}^2$. It may be
$I_s=\varnothing$, if $n_s=0$. Two indexes $\a,\b\in\{1,...,n\}$
are said to be equivalent, and we write $\a\sim\b$, if they belong
to the same set $I_j$.
 Two indexes $\a,\b$ are said to be orthogonal, and
we write $\a\perp\b$, if $\a\in I_j$ and $\b\in I_k$ with $k\neq j$.

Let
\begin{equation}\label{dominio}
\begin{aligned}
 &\Omega:= \big\{ (z,z^{n+1})\in \C^{n+1}  \, : \, \Im(z^{n+1})>
p(z,\zz)\big\},
\\
 & M_0 :=b\Omega := \big\{(z,z^{n+1}) \in \C^{n+1} \, : \, \Im(z^{n+1})=
p(z,\zz) \big\}, \qquad\text{and}
\\
  &M  : = \Big\{(z,z^{n+1})\in M_0 :\prod_{j=1}^{s-1}\abs{z_j}\neq 0 \Big\}.
\end{aligned}
\end{equation}
We will see that $M$  is the strictly pseudoconvex part of the surface $M_0$.
The unbounded open set $\Omega $ is biholomorphically equivalent to
the bounded generalized ellipsoid
\[
 E:=\left\{(w, w^{n+1})\in \C^{n+1}: \sum_{j=1}^s \abs{w_j}^{2m_j} + \abs{w^{n+1}}^2<1\right\}
\]
via the map
\begin{equation*}
 \Omega\ni (z, z^{n+1})\longmapsto
\Big(\frac{2^{1/m_1}z_1}{(i+z^{n+1})^{1/m_1}},
,\dots,\frac{2^{1/m_{s-1}}z_{s-1}}{(i+z^{n+1})^{1/m_{s-1}}},  \frac{2
z_s}{i+z^{n+1} }, \frac{i-z^{n+1}}{i+z^{n+1}}\Big) \in E.
\end{equation*}
In the rest of the paper, we will work on the unbounded model, where the
computations are easier.

\begin{remark}
In \eqref{ipo} we require  $m_j\ge 2$ and $n_j\ge 2$ for $j=1, \dots,
s-1$.
The case    $n_j=1$ for all $j$ is  discussed  in the recent paper
\cite{Landucci}.
Assumption \eqref{ipo} ensures that all local CR mappings extend to
global biholomorphic mappings, see \cite{DS}. If \eqref{ipo} is
violated, this is in general not  true.
\end{remark}

Let us consider  the  following biholomorphic  mappings.
The mapping $I:\Omega\to\Omega$
\begin{equation}
 I(z_1,\dots, z_{s-1}, z_s, z^{n+1})   =  \Big(\frac{z_1}{(z^{n+1})^{1/m_1}},
\dots,
\frac{z_{s-1}}{(z^{n+1})^{1/m_{s-1}}} ,\frac{
z_s}{z^{n+1}}, -\frac{1}{z^{n+1}}\Big) \label{mod2}
\end{equation}
is the inversion. For any $r>0$, the mappings $\delta_r:\Omega\to\Omega$,
\begin{equation}
   \d_r(z_1,\dots, z_{s-1},  z_s, z^{n+1})   =
\big(r^{1/m_1}z_1,\dots,r^{1/m_{s-1}}z_{s-1} , rz_s, r^2
z^{n+1}\big)\label{mod3}
\end{equation}
form a one-parameter group of dilatations. Finally, consider the mappings $\psi$
 of the form
\begin{equation}
 \psi(z,z^{n+1})=\Big( B_1 z_{\sigma(1)},   \dots, B_{s-1}z_{\sigma(s-1)},
B_sz_s + b_s,
  b^{n+1}+ z^{n+1}+2i  ( B_sz_s \cdot  \bar b_s)\Big), \label{mod4}
\end{equation}
where
$\sigma$ is a permutation  of $\{1, \dots, s-1\}$ such that
              $m_{\s(j)}= m_{j}$ and  $n_{\s(j)}= n_{j}$ for any $j=1, \dots,
              s-1$, $B_j\in U(n_j)$ are unitary matrices, $b_s\in \C^{n_s}$,
and $b^{n+1}= t_0+ i\abs{b_s}^2\in\C$ for some $t_0\in\R$.

For $a=(a_s,a^{n+1})\in C^{n_s}\times\C$ with $a^{n+1} = t_0+i |a_s|^2$ for
some $t_0\in\R$, let $\phi_a$ be the mapping
\begin{equation}
 \phi_a(z_1,\dots,z_{s-1}, z_s, z^{n+1})   = (z_1,\dots, z_{s-1}, z_s+a_s,
z^{n+1}+a^{n+1}+2iz_s\cdot\bar a_s) \label{mod1}.
\end{equation}
The mapping $\phi_a$ is a
particular case of  \eqref{mod4}.

A composition of the mappings \eqref{mod2}--\eqref{mod4}
extends to a CR mapping of $M_0$, possibly off one point. Our main theorem
states that any local $CR$ mapping of $M$ is such a
composition.

\begin{theorem}
 \label{principale}Let $f:N\to N'$ be a CR mapping between connected open
subsets of $M$. Then,
for a suitable choice of $\psi$ as in \eqref{mod4}, $r>0$, and $a
=(a_s,t_0+i|a_s|^2)\in \C^{n_s}\times\C$ we have
\begin{equation}
 f= \psi\circ \d_r \circ J\circ \phi_a ,
\end{equation}
where either $J= I$ as in in \eqref{mod2} or $J$ is  the
identity map.
\end{theorem}

\begin{proof}[Scheme of the proof of Theorem \ref{principale}] $\,$

\emph{Step 1. } For a suitable contact form $\theta$ on $M$, we
show that the CR factor $\l_f$  of a CR function $f$ between open,
connected subsets of $M$ either is a constant or
 it has the form
\begin{equation}\label{eccolaone}
\l_f =  k^{-2}
\abs{z^{n+1}+ a^{n+1}+2iz_s\cdot \bar a_s}^{-2}, \end{equation}
 for some $k>0$,  $a_s\in\C^{n_s}$ and
$a^{n+1} = t_0+i \abs{a_s}^2\in \C$, where $a_s=0$ if $n_s=0$.
 This is proved in Theorem \ref{oneone}. The proof requires the study of the
overdetermined system   in  \eqref{ochetta} and
\eqref{pablo}. To solve the system, we exploit the structure of
the Chern-invariant
cone bundle $\H$. This is carried out in Subsection \ref{imparo}.

\emph{Step 2.}  Once the form of the CR factor $\lambda_f$ is known,
 we consider the mappings $\phi_a$
in
\eqref{mod1} and $\delta_r$ in \eqref{mod3} with $r=1/k$. Elementary
computations on the CR factors give:
\begin{equation}\label{confo}
 \l_{\phi_a}(z)  = 1,\qquad \l_{I}(z)=  |z^{n+1}|^{-2}, \quad
\text{and}\quad \l_{\d_r} = r^2.\end{equation}
Let $G:=\delta_{1/k} \circ I\circ\phi_a$ and define the mapping $\psi$
via the identity $f=
\psi\circ G$.
By \eqref{confo}, $\psi$ satisfies $\psi^*\theta=\theta$.
Indeed, the CR factor $\l_G$ of $G$ is
\begin{equation*}
 \l_{G}(z) = \l_{\d_{1/k}}(
I(\phi_a(z))\l_I(\phi_a(z))\l_{\phi_a}(z)=  k^{-2}\big|z^{n+1}+
a^{n+1}+ 2iz_s\cdot\bar a_s\big|^{-2},
\end{equation*}
and therefore
\begin{equation*}
 \l_f(z)=\l_\psi(G(z))\l_G(z)=\l_\psi(G(z)) k^{-2} \big|z^{n+1}+
a^{n+1}+ 2iz_s\cdot\bar a_s\big|^{-2}.
\end{equation*}
Thus, by the form of $\l_f $ in \eqref{eccolaone} we deduce that the CR
factor of $\psi$ is  $\l_\psi=1$.
In Subsection \ref{leviso} we show that all such mappings, that we call
\emph{Levi-isometric}, have the form \eqref{mod4}.

This concludes the proof of the classification result.
\end{proof}

\section{Pseudohermitian and Chern--Moser invariants in generalized ellipsoids}
\label{invarianti}

\setcounter{equation}{0}
\subsection{Computation of the pseudohermitian invariants}\label{cierre}
 Fix on $M_0$ the pseudohermitian structure $\vartheta =
i(\p F- \bar \p F)$, where $F(z,\bar z, z^{n+1},\bar z^{n+1})=
\mathrm{Im}(z^{n+1}) - p(z,\bar z)$ is a defining function for
$M$:
\begin{equation}
 \label{conto}
 \vartheta : = \frac{dz^{n+1} +d\zz^{n+1}}{2}-i(p_\a dz^\a- p_\aa d\zz^\a)
           = dt -i(p_\a dz^\a- p_\aa d\zz^\a).
\end{equation}
We use the notation  $p_\a  =\frac{\p p}{\p z^\a}$ and we let $t=\Re(z^{n+1})$.
  On $M_0$ we fix the coordinates
$(z,t)\in\C^n\times\R$, i.e., we identify $(z,t)\in\C\times\R$ with
$\big(z,t+ip(z,\bar z)\big) \in M_0$.
Fix the  holomorphic frame
\begin{equation}\label{zaza}
    Z_\a  =
         \p_\a+ip_\a\p_t\qquad \text{for
$\a=1,...,n,$}
                     \end{equation}
and let $Z_\aa = \overline{ Z}_\a$. We clearly have $\vartheta(Z_\a)=\vartheta(Z_\aa)=0$ for any $\a=1,...,n$.

 The Levi form on $M$ is the $2$-form $L =-i d\vartheta$.
From the identities
\begin{align}
 \label{ACCA}
   & h_{\a\bb} :=   L(Z_\a, Z_\bb)=-id\vartheta (Z_\a,Z_\bb) =
i\vartheta([Z_\a, Z_\bb]) \quad\text{and}
   \\ \label{kom1}
& [Z_\a, Z_\bb]= -2ip_{\a\bb}\p_t\quad   \qquad \text{for
$\a,\b=1,\dots,n,$}
\end{align}
we obtain
 \begin{equation}
 \label{hab}
   h_{\a\bb}  =
  \begin{cases}
  \displaystyle
    2m_j |z_j|^{2(m_j-1)} \Big( \d_{\a\b}
    +(m_j-1)\frac{\bar z^\a z^\b}{|z_j|^2}\Big)
     & \text{if }\a\sim \b\in I_j, \\
    0 & \text{if } \a\perp\b.
  \end{cases}
\end{equation}
The inverse  matrix $h^ {\la\bb}$ has the form
\begin{equation}
 \label{hab1}
   h^ {\la\bb} = \begin{cases}
  \displaystyle
    \frac{1}{2m_j |z_j|^{2(m_j-1)}}
     \Big( \d_{\la\b}
    -\frac{m_j-1}{m_j}
    \frac{z^\l\bar  z^\b}{|z_j|^2}\Big)
     & \text{if }\la\sim \b\in I_j, \\
    0 & \text{if } \l\perp\b.
  \end{cases}
\end{equation}

The surface $M$ defined in \eqref{dominio}
 is strictly pseudoconvex,
 because $\det h_{\a\bb}>0$ on $ M$.
The characteristic vector field  of $\vartheta$ is
$     T = \p_{n+1} + \bar \p_{n+1} = \frac{\p}{\p t}.
$

For   $j=1,...,s$,  let us introduce the holomorphic vector fields
\begin{equation}\label{cirio}
  E_j = \frac{1}{m_j} \sum_{\a \in I_j} z^\a Z_\a.
\end{equation}
We say that $E_j$ is a  vector field ``of
radial type''.
A short computation based on \eqref{hab} shows that, for any  $\a\in\{1,\dots ,n\}$, $j,k\in\{1,\dots,s\}$,
\begin{equation}\begin{aligned}
 L(Z_\a, \overline E_j) &= 2m_j |z_j|^{2(m_j-1)}\zz^\a
\qquad \text{ if $\a\in I_j$,\quad  and}
 \\  L(E_j, \overline E_k)& = 2 \abs{ z_j}^{2m_j}\d_{jk} \qquad\text{for all $j,k=1, \dots, s$.}
\label{apparato}
\end{aligned}\end{equation}
Here, $\delta_{jk}$ is the Kronecker's symbol. The vector fields $E_1,...,E_s$
 span an $s$-dimensional subbundle $\mathcal E\subset T^{1,0}M$.
   If $n_s=0$, we have no vector field $E_s$ and $\E$ has
dimension $s-1$.  We denote
by $\mathcal
E^\perp$ the orthogonal complement of $\mathcal E$ in $\T M$ with
respect to the Levi form.

Let $Q:\T  M\to \E^\perp$ be the projection
\begin{equation}\label{q}
  Q(Z)= Z-\sum_{j=1}^s \frac{L(Z, \EE_j)}{L(E_j, \EE_j)}E_j.
\end{equation}
  If $n_s=0$, in the sum the index $j$ ranges from $0$ to $s-1$.
In particular,  for
any $j\in\{1,\dots,s\}$ and  $\a\in I_j$, let  $W_\a $ be the holomorphic vector field
\[
   W_\a : =  Q(Z_\a)
          = \sum_{\b\in I_j}Q_\a^\b Z_\b.
\]
By \eqref{apparato}, \eqref{q}, and \eqref{cirio}, we deduce that the
coefficients $Q_\a^\b$ are
\begin{equation}
 \label{carie} Q_\a^\b = \d_\a^\b -\frac{\bar z^\a z^\b}{\abs{z_j}^2}.
\end{equation}
Let us introduce the hermitian form $Q_\flat$ on $T^{1,0}M$ associated with $Q$:
\begin{equation}
  \label{FLAT}
 Q_\flat(U, \VV) : = L(Q(U), \VV)=L(Q(U), \overline{Q(V)})\qquad
\text{for any }  U, V\in \T M.
\end{equation}
Letting $Q_{\a\bb} := Q_\flat (Z_\a,Z_\bb)$, we have for $\a\sim\b\in I_j$,
\begin{equation}
  \label{sette}
 Q_{\a\bb} = Q_\a^\g h_{\g\bb}=
       2m_j\abs{z_j}^{2(m_j-1)}\Big(   \d_{\a\b} - \frac{\zz^\a z^\b}
 {\abs{z_j}^2 }\Big)
 =2m_j\abs{z_j}^{2(m_j-1)}Q_\a^\b,
\end{equation}
whereas $Q_{\a\bb}=0$  if $\a\perp\b$.

Finally observe  that \eqref{zaza} and \eqref{carie}  give
\begin{equation}
 \label{sfer}
  W_\a   =\p_\a -  \sum_{\b\in I_j} \frac{\bar z^\a z^\b }{|z_j|^2}
  \p_\b\qquad \text{if $\a\in I_j$.}
\end{equation}
Thus $W_\a p=0$ for $\a=1,...,n$, i.e., $W_\a$ is a holomorphic
vector field in $\C^n$ which is tangent to the hypersurfaces of
$\C^n$ given by  $p(z,\bar z) = $  constant.

The bundle $\mathcal  E^\perp$ is $n-s$ dimensional  (in fact,
($n-s+1)$--dimensional if $n_s=0$) and it is
generated by the vector fields  $W_\a$ with  $\alpha=1,...,n$.
Let $\mathcal W_j\subset \T M$ be the subbundle spanned by the
vector fields $W_\a$ with $\a\in I_j$, where $\W_s=(0)$, if $n_s=0$ or $n_s=1$.
Then we have
$  \mathcal E^\perp = \mathcal W _1 \oplus\cdots\oplus \mathcal{W}_{s-1}\oplus\mathcal
  W_s.
$
Finally,  we have $\W_j\perp\W_k$ if $j,k\in \{1,\dots,s\}$ are different.
Indeed, for any  $\a\in I_j$ and $ \b\in I_k$, we have
$
 L(W_\a, W_\bb)= Q_\a^\g Q_\bb^\ss h_{\gamma\ss} =0,
$
because $h_{\g\ss}=0$ if $\g\perp\s$.
Therefore the decomposition
\begin{equation}\label{caimano}
  \T M =\W_1\oplus\cdots\oplus\W_{s-1}\oplus\W_s\oplus\E
\end{equation}
is orthogonal.
Observe also that $\W_s\oplus \E = \mathrm{span}  \{ E_1, \dots, E_{s-1}, Z_\a : \a\in I_s\}$.

\begin{proposition}
 \label{qam}
The hermitian form $Q_\flat$  introduced in \eqref{FLAT} satisfies:
\begin{subequations}
\begin{align}
 \label{Q1}
 &
  Q_\flat(E,\bar Z)  = 0\qquad \text{for all $E\in\E$ and
  $Z\in \T M$;}
  \\ \label{Q2}
  &
  Q_\flat(V, \bar W) =
   L(V,\bar W)\qquad \text{for all $j\in\{1, \dots, s\} $ and $  V,W\in \mathcal W_j$;}
 \\ \label{Q3}
  &
  Q_\flat(Z,\bar Z) \geq 0\quad  \text{ for all $Z\in \T M$\quad  and}\quad
  \mathcal E = \big\{Z\in \T M : Q_\flat(Z,\bar Z) = 0 \big\}.
\end{align}
\end{subequations}
\end{proposition}

\begin{proof}
Let $E\in\E$ and $Z\in \T M$. Then
$
 Q_\flat(E, \ZZ)=L(Q(E), \ZZ)=0,
$
because $Q(E)=0$. This proves the first line.
To prove the second line, just observe that
$
 Q_\flat (V, \WW)= L(Q(V), \WW)= L(V, \WW),
$
because $Q(V)=V$. The third line \eqref{Q3} follows from \eqref{FLAT}, letting $U=V$
and from strict pseudoconvexity.
\end{proof}

Let $\D$ be the Tanaka-Webster connection of  $  (M, \theta)$.
We refer to \cite{Tanaka75,W77} and \cite{DT} for the relevant facts concerning $\D$.
The curvature operator of $\D$ is
$    R(Z_\la,Z_{\bar\mu}) Z_\a = \D_{Z_\la} \D_{Z_{\bar\mu}} Z_\a-
   \D_{Z_{\bar\mu}} \D_{Z_\la} Z_\a -\D_{[Z_\la, Z_{\bar\mu}]} Z_\a.
$
The   curvature  tensor have components
$   R_{\a\bar\b\la\bar\mu} =
         L ( R(Z_\la,Z_{\bar\mu})Z_\a,Z_{\bar\b}).
$ It
enjoys the symmetries
\begin{equation}
 \label{sym}
   R_{\a\bar\b\ga\bar\mu}=R_{\ga\bar\b\a\bar\mu}\quad\text{and}\quad
  R_{\a\bar\b\ga\bar\mu} =  \bar R_{\b\bar\a\mu\bar\ga},
\end{equation}
see \cite[Section 1.4]{DT}.
The pseudohermitian  Ricci tensor is defined by
$
   R_{\a\bar\mu} = R_{\a\,\,\,\la\bar\mu}^{\,\,\,\la},
$
and  the scalar curvature is $ R =
h^{\a\bar\mu} R_{\a\bar\mu}$.
 Finally, the pseudohermitian torsion of $\D$ is defined by
$
 \t (Z_\b) = \D_T Z_\b - \D _{Z_\b }T - [T, Z_\b] = :A_\b^\aa Z_\aa,
$
where
$
   A_{\a\b} :=    L(\tau(Z_\a),Z_\b)
$  satisfies $A_{\a\b}=A_{\b\a}$, as proved in \cite{W77}.


Now we study the curvature tensors on the hypersurface
$M. $  Associated with the
decomposition $\C T   M = T^{1,0}   M \oplus T^{0,1}  M \oplus \C T$, we
have the projections $\pi_+:\C T   M\to T^{1,0}  M$ and
$\pi_-:\C TM \to T^{0,1}M$. By definition, for $U,V $ holomorphic vector fields,
we have $\D_U \bar V : = \pi_-([U,\bar V])$. In our case, from
\eqref{kom1} we find
\begin{equation}
 \label{D1}
   \D _{Z_\a} Z_{\bar\b} 
   =0\qquad \text{for all $\a,\b=1,...,n$.}
\end{equation}

If $U, V $ are holomorphic vector fields, then $\D_U V  $ is  defined by
 $L(\D_U V, \bar W) = U
L(V,\bar W) - L(V, \D_U\bar W)$ for all $W\in T^{1,0}M$. Thus, we obtain
\begin{equation}
 \label{lixor}
      \D_{Z_\la} Z_\a = (h^{\s\bar\b} \p_\la h_{\a\bar\b}) Z_\s \qquad\text{for all  $\a,\la=1,...,n.$}
\end{equation}
Since  $h_{\a\bar\b} = 2 p_{\a\bar\b}$, we deduce from
\eqref{lixor} that
\begin{equation}
 \label{orthus}
   \D_{Z_\la} Z_\a = 0 \qquad\text{for all $\a,\la\in 1, \dots,n$ with $\a\perp\lambda$.}
\end{equation}
The Tanaka--Webster connection satisfies $\D T=0$. Moreover, in our
case we have $ \D_T Z_\a = \pi_+([T,Z_\a]) = 0$. By \eqref{kom1},
we deduce that
\begin{equation}
 \label{lalix}
       \D_{[Z_\la,Z_{\bar\mu}]} Z_\a = 0 \qquad \text{for all $ \a, \l,\m\in \{1,\dots,n\}$.}
\end{equation}
By \eqref{D1} and \eqref{lalix}, the curvature operator
reduces to
$
   R(Z_\la,Z_{\bar\mu}) Z_\a
   =-
   \D_{Z_{\bar\mu}} \D_{Z_\la} Z_\a$,
and taking into account \eqref{lixor}, we get the Riemann and Ricci  tensors
\begin{align}
\label{rullo}
  & R_{\a\bar\b\la\bar\mu}   = - \p_{\bar\mu} (h^{\s\bar\ga} \p_\la
        h_{\a\bar\ga})  h_{\s\bar\b},
  \\
&  R_{\a\bar\mu}  = R_{\a\,\,\,\la\bar\mu}^{\,\,\,\la}
    = - \p_{\bar\mu} (h^{\la \bar\ga} \p_\la
        h_{\a\bar\ga}).
\label{rallo} \end{align}
Finally, since   $\D_{Z_\a} T = \D_T Z_\a = [Z_\a, T]=0$, the torsion vanishes identically,
\begin{equation}
 \label{tower}
  A_{\a\b} = 0.
\end{equation}

\begin{proposition}
      \label{mirox}
 The curvature tensor  of $\D$ of $(M,\vartheta)$
 has the form
 \begin{equation}
  \label{detto}
  R_{\a\bb\l\mm}= - \frac{m_j-1}{2m_j}|z_j|^{-2m_j} \big\{
       Q_{\l\mm}Q_{\a\bb}+ Q_{\a\mm}Q_{\l\bb}\big\}
       \qquad\text{if }
       \a,\b,\l,\m\in I_j,
\end{equation}
for some $j\in\{1, \dots, s\} $, and $R_{\a\bb\l\mm}=0$ if two of the indeces
$\a,\b,\lambda,\mu$ are orthogonal.
Moreover, $R(U, \bar V, Z,
\WW)=0$ as soon as one of the vector fields $U,V,Z,W\in T^{1,0}M$
belongs to  $\mathcal E\oplus\mathcal{W}_s  $.
The pseudohermitian Ricci tensor
 has the form
\begin{equation}
\label{ture}
  R_{\l\mm}
   =-\frac{n_j(m_j-1)}{2m_j}|z_j|^{-2m_j} Q_{\l\mm}
   \qquad \text{if }\l,\m \in I_j,
\end{equation}
and $R_{\l\mm}=0$ if $\l\perp\mu$. Moreover, $\mathrm{Ric}(U,\bar
V)=0$ as soon as one of the vector fields $U,V\in T^{1,0}M$ belongs to
$\mathcal E\oplus\W_s$.
Finally, the scalar curvature is
\begin{equation}
 \label{coccobill}
 R= -\sum_{j=1}^{s-1}\frac{n_j(n_j-1)(m_j-1)}{2m_j}\abs{z_j}^{-2m_j}.
\end{equation} \end{proposition}

\begin{proof}
We start from the formula \eqref{rullo}. The components of the Levi form are given
in \eqref{hab} (and \eqref{hab1}). Note that
\[
  \begin{array}{lll}
    \a\perp\la  & \Rightarrow & \p_\la h_{\a\gg} = 0,
      \\
     \a\perp \mu & \Rightarrow  & \p_{\bar\mu} (h^{\s\bar\ga} \p_\la
        h_{\a\bar\ga})=0,
        \\
     \a\perp \b & \Rightarrow &
        \p_{\bar\mu} (h^{\s\bar\ga} \p_\la
        h_{\a\bar\ga})  h_{\s\bar\b} =0.
  \end{array}
\]
Using the symmetries \eqref{sym}, we conclude that $
R_{\a\bb\l\mm}=0$ as soon as there are two orthogonal indexes.

Assume that $\a,\b,\l,\m$ are in $I_j$. From \eqref{hab} we
get
\[
 \p _\l h_{\a\gg}
 = 2m_j(m_j-1) |z_j|^{2(m_j-2)}
   \Big\{ \zz^\l \d_{\a\g}+\zz^\a \d_{\l\g}
   +(m_j-2) \frac{\zz^\a z^\g \zz^\l}{ |z_j|^2} \Big\},
\]
and thus
\[  h^{\s\gg} \p_\l h_{\a\gg} = \frac{m_j-1}{|z_j|^2}
   \Big\{ \zz^\l \d_{\a\s}+\zz^\a \d_{\l\s}
    -\frac{\zz^\a z^\s \zz^\l}{ |z_j|^2} \Big\}.
\]
After a short computation based on  \eqref{carie} and \eqref{sette}, we find
\begin{equation}
  \label{Pump}\begin{aligned}
-R_\a{}^\s{}_{\l\mm}&=   \p_{\mm}\big(h^{\s\gg} \p_\l h_{\a\gg} \big) =
  \frac{ m_j-1}{\abs{z_j}^{2}}
     \big\{ Q_\a^\s  Q_{\la }^\m+ Q_\la^\s Q_{\a}^\m\big\}
     \\&
     = \frac{m_j-1}{2m_j}\abs{z_j}^{-2m_j}
   \big\{ Q_\a^\s  Q_{\la\mm}+ Q_\la^\s Q_{\a\mm}\big\},
              \end{aligned}
\end{equation}
and contracting with $h_{\s\bb}$, we get \eqref{detto}.

Next we show that $R(Z, \WW, U, \VV)=0$ if $Z\in \W_s\oplus\E$. If
$Z\in \E$, this follows trivially from the first line of
\eqref{Q1}. If  $Z\in \W_s$, then $R(Z, Z_\bb, Z_\l, Z_\mm)=0$ if
at least one of the indexes $\b,\l,\m$ does not belong to $I_s$.
If all $\b,\l,\m\in I_s$, then, by \eqref{detto}, $R(Z, Z_\bb,
Z_\l, Z_\mm)=0$, because $m_j-1=0$ if $j=s$.

Identities \eqref{ture} and \eqref{coccobill} follow upon contracting
indexes in \eqref{detto}, $R_{\l\mm}=  R_\a{}^\a{}_{\l\mm}$. Recall that
  by \eqref{carie}
  and \eqref{sette}
   we have $\sum_{\a\in I_j}Q_\a^\a= n_j-1$ and $\sum_{\s\in I_j}Q_\a^\s
Q_{\s\mm}= Q_{\a\mm}$, if $\a,\m\in I_j$.
\end{proof}

\begin{remark}
 \label{miros}
 Let $V\in \mathcal E^\perp$ be a vector  such that
 $V=V_1+...+V_{s-1}$  with $V_j\in  \mathcal W_j$. The pseudohermitian sectional
 curvature of $(M,\vartheta)$ along $V \neq 0$ is
 \begin{equation}
  \label{sec}
     k(V):  = \frac{R(V,\bar V,V,\bar V)}{|V|^4} = - \frac{1}{|V|^4}
     \sum_{j=1}^{s-1} \frac{m_j-1}{m_j} |z_j|^{-2m_j} |V_j|^4,
\end{equation}
where $\abs{V}:= L(V, \VV)^{1/2}$ denotes the Levi-lenght of $V$.
This formula follows from \eqref{detto} and  \eqref{Q2}. Notice, in particular, that, since   $m_j>1$ for all
$j\le s-1$,  then  $k(V)\neq 0$ for any $V\in \mathcal \W_1\oplus\dots\oplus\W_{s-1}$ with $V\neq
0$.
\end{remark}

\subsection{Chern invariant cones in generalized ellipsoids}
\label{imparo}
We describe the structure of the cones   $\H_P$ introduced in \eqref{accadi}.
Here and hereafter, let $\abs{U}:= L(U, \UU)^{1/2}$ denote the
Levi-length of $U\in \T M$.

\begin{proposition}
  \label{pss}
 Let $M\subset\C^{n+1}$ be
the surface defined in \eqref{dominio}.
 Then  \begin{equation}
 \label{Cinv}
\H = \mathcal W_1\cup ...\cup \W_{s-1}\cup\big( \mathcal W_s \oplus   \E\big).
\end{equation}
\end{proposition}

\begin{proof}
We prove that $\E \oplus\W_s \subset \H$. In fact, if $U\in \E\oplus\W_s$,  then $ R(U,
\VV,  Z, \WW) = 0$ for all $V,Z,W\in\T M$, by Proposition
\ref{mirox}.

In order to show that $\mathcal W_j\subset\H$ for all $j=1,\dots,s-1$,  let
 $U\in \mathcal W _j$ and take $V,Z,W\in  {T^{1,0}M}$ 
such that $L(U, \VV)= L(U, \WW)= L(Z, \VV)= L(Z, \WW)=0$. Observe that, writing $V= V_j+V_j^\perp$, where $V_j\in\W_j$ is the projection of $V$ onto $\W_j$ and $V_j^\perp= V- V_j= E+ \sum_{k\neq j}V_k$, for some $E\in\E$, by \eqref{caimano},
 we have
\begin{equation}
 \label{jk} L(U, \VV)= L(U, \VV_j) \quad\text{and}\quad L(U, \WW)= L(U, \WW_j).
\end{equation}
Here, we made for $W$ the same decomposition as for $V$. Observe also that
$Q_\flat(U,\VV ) = Q_\flat(U, \VV_j),$ and   $Q_\flat(U,\WW ) = Q_\flat(U, \WW_j)$, by \eqref{FLAT}.
Then
we have by \eqref{detto}
\[
\begin{split}
 R(U, \VV,  Z, \WW)
 & = - \frac{m_j-1}{2m_j} |z_j|^{-2m_j} \big\{
       {Q_\flat} (U,\VV_j)Q_\flat(Z,\WW)+ Q_\flat(U,\WW_j) Q_\flat(Z,\VV)
\big\}
       \\
       &
       = - \frac{m_j-1}{2m_j} |z_j|^{-2m_j} \big\{
       L(U,\VV_j)Q_\flat(Z,\WW)+ L(U,\WW_j) Q_\flat(Z,\VV) \big\}
       =0,
\end{split}
\]
by \eqref{jk} and by the definition of $\H$.

Now we show that any vector field $U = X+Y$ with $X\in \E\oplus\W_s$ and
$Y\in \W_1\oplus\cdots\oplus\W_{s-1} $ such that $|X|\neq 0$ and $|Y|\neq0$ does not
belong to $\H$. We assume without loss of generality that $|Y|=1$.
We choose $Z = X+Y$ and $V=W=X-k Y$ with $k =|X|^2\neq 0$, in such
a way that $L(U,\VV) = L(X+Y,\bar X- k \bar Y) = 0$. By
Proposition \ref{mirox},
$
 R(U,\VV,Z,\WW) = k^2 R(Y,\bar Y, Y,\bar Y)= k^2|Y|^4 k(Y) \neq 0,
$
thanks to  \eqref{sec}.

Finally, we prove that if $V\in (\E\oplus\W_s)^\perp\setminus (\mathcal
W_1\cup \cdots \cup \mathcal W_{s-1})$ then $V\notin \H$. Let
$V=V_1+\cdots +V_s$ with $V_j\in \mathcal W_j$ and assume without loss
of generality that $|V_1|\neq0$ and $|V_2|=1$. Let $W=V_1-\kappa
V_2$ where $\kappa = |V_1|^2$ in such a way that $L(V,\WW) = 0$.
By Proposition \ref{mirox} and Remark \ref{miros}, we have
$
   R(V,\WW,V,\WW) = R(V_1,\bar V_1,V_1,\bar V_1)+\
   \kappa^2 R(V_2,\bar V_2,V_2,\bar V_2)\neq 0,
$
as claimed. \end{proof}

Let $N, N'\subset   M$ be connected open sets and let $f:N\to N'$ be a
Cauchy--Riemann diffeomorphism. By
\eqref{accado} it must be $f_*(\H_P)=\H_{f(P)}$. Then there are two
cases:

\textbf{(A)}  $f_*\mathcal (\E\oplus\W_s) = \mathcal E \oplus\W_s $;

\textbf{(B)} there exists $j=1,..., s-1$ such that $f_*(\E\oplus\W_s) =
\mathcal W_j $.

\noindent Here and hereafter, with slight
 abuse of notation let  $f_*\mathcal (\E\oplus\W_s) = \mathcal E \oplus\W_s $
stand for
$f_*\mathcal (\E\oplus\W_s)_P  = (\mathcal E \oplus\W_s)_{f(P)} $, for all $P\in N$.
\medskip

\noindent  Case (B) may occur only if $\mathrm{dim}(\E\oplus\W_s)
= \mathrm{dim}(\mathcal W_j)$. Actually, the Case B cannot occur
at all, as the following theorem states.

\begin{theorem}
 \label{bad}
 Let $N, N'\subset M$ be   open sets.
 A CR diffeomorphism $f:N\to N'\subset M$ satisfies
 $f_*(\E\oplus\W_s)= \E\oplus\W_s$.    In particular, there exists a permutation $\sigma$
 of $\{1,...,s-1\}$ such that $f_* (\W_j ) \subset  \W_{\sigma(j)}$ for
all
 $j=1,...,s-1$.
 \end{theorem}

We   prove Theorem \ref{bad} in Section \ref{brutta}.
In the following proposition,
 we prove   that, in case (A), CR diffeomorphisms   preserve the Ricci
Tanaka--Webster curvature of $\theta$.
Diffeomorphisms that preserve the Ricci curvature of a Riemmanian metric are
rather studied in Riemannian geometry, see
\cite{Kuhnel,Osgood}. It could  be  interesting to see whether such
  mappings enjoy any significant geometric property in the CR setting.

\begin{proposition}
  \label{ricc}
 Let $N, N'\subset M$ be  open sets. A CR
  diffeomorphism
 $f:N\to N'$ such that $f_*(\E\oplus\W_s)= \E\oplus\W_s$
  preserves the Ricci curvature of
  $\vartheta$, i.e.
   $\Ric(f_*
   Z, f_*\bar W)=\Ric(Z, \WW)$ for all $Z, W\in \T N.$
\end{proposition}
\begin{proof}
We divide the proof into two steps.

\it Step 1. \rm We claim that $R = \la R\circ f$, where $R$ is the
scalar curvature of  $\vartheta$ and $\la>0$ is the
CR
factor of $f$, i.e., $f^*\vartheta = \la\vartheta$.

\medskip

Let $Z\in\mathcal E\oplus\W_s$ with $|Z|=1$.
Proposition \ref{mirox} and   Webster formula \eqref{web} yield
\begin{equation}
 \label{omega}
 S(Z,\bar Z, Z,\bar Z) = \frac{2R}{(n+1)(n+2)}.
\end{equation}
Since $f_*Z\in\E\oplus\W_s$, we analogously have
\begin{equation}
 \label{voromega}
 S(f_*Z,f_*\ZZ, f_*Z,f_*\ZZ )
   = \frac{2R\circ f}{(n+1)(n+2)}|f_* Z|^4
   =   \frac{2\la^2 R\circ f}{(n+1)(n+2)},
\end{equation}
where we also used $|f_* Z|^2 = \la |Z|^2=\la$.  Thus,
by the relative CR invariance \eqref{chevi} we have
\[
S_\theta(f_*Z, f_*\ZZ, f_*Z, f_*\ZZ) =
S_{f^*\vartheta}(Z, \ZZ, Z, \ZZ)=    S_{\lambda \vartheta}(Z, \ZZ, Z, \ZZ) = \la
S_\theta(Z, \ZZ, Z, \ZZ) ,
\] Here,
$S_\theta, S_{f^*\vartheta}$ and $S_{\lambda\vartheta}$ denote  the
Chern
tensors relative to $\theta, f^*\vartheta$, and $\lambda\vartheta$,
respectively. Comparing \eqref{omega} and \eqref{voromega} we conclude  \it Step~1\rm.

\medskip

\it Step 2. \rm We claim that  $\mathrm{Ric}(f_* V,f_*\bar W) =
\mathrm{Ric}(V,\bar W)$ for all $V,W\in (\mathcal E \oplus W_s)^\perp$.

\medskip

Take vector fields $V,W\in (\E \oplus W_s)^\perp$. Let also $Z\in \E \oplus W_s$ be such that $|Z|=1$. All terms in the Chern
tensor containing curvature tensors along $Z$ or  terms of the form
$L(V,\bar Z)  $ and $ L(Z,\WW)$ vanish. Thus
\begin{equation}
 \label{silix}
   S(V,\bar W, Z, \bar Z) = -\frac{1}{n+2} \left\{
   \mathrm{Ric}(V,\bar W)
      - \frac{R}{(n+1)} L(V,\bar W) \right\}.
\end{equation}
Since $f_* Z \in \E \oplus W_s$ and $f_* W ,f_*V\in (\mathcal E \oplus W_s)^\perp$, we
analogously have
\begin{equation}
 \label{selene}
\begin{split}
    S(f_* V,f_*  \bar W, f_* Z,  f_* \bar Z)
     & = -\frac{\la }{n+2}\left\{
      \mathrm{Ric}( f_* V, f_*\bar W)-
      \frac{\la R\circ f }{n+1} L(V, \bar W)\right\}.
\end{split}
\end{equation}
Recall that
$  S(f_* V,f_*  \bar W, f_* Z,  f_* \bar Z) =   S_{\lambda\vartheta}(V,\bar W,
Z,
  \bar Z)=\la S(V,\bar W, Z,  \bar Z),
$ by the CR invariance \eqref{chevi}. Thus  the \emph{Step 2} can be
accomplished on
comparing  \eqref{silix}, \eqref{selene} and using the    \emph{Step 1}.

 The proof is finished, because $\Ric(Z, \WW)=0$
 for all $Z\in\E\oplus\W_s$ and $W\in\T  M$.
\end{proof}

\begin{proposition}
 \label{radioso}
Let $N, N'\subset M$ be  open sets and let $f:N\to N' $ be a
CR diffeomorphism such that $f_*(\E\oplus \W_s) =
\E\oplus \W_s$. Then the CR factor $\la$ of $f$ satisfies
$W_\a \l =0$ for any  $\alpha\in I_1\cup\cdots \cup I_{s-1}$.
  \end{proposition}

\begin{proof}
Let $\alpha\in I_1\cup\cdots\cup I_{s-1}$ and fix $W = W_\alpha$. We
first observe that
\begin{equation*}
\begin{aligned}
 d\theta(f_*T, f_*\WW)&= f^*(d\theta)(T, \WW) = d(f^*\theta)(T, \WW) =
\big((d\lambda)\wedge\theta + \lambda d\theta\big)
(T, \WW)=-\WW\lambda.
\end{aligned}
\end{equation*}
Therefore, it suffices to show that $d\vartheta(f_*T, f_*\WW)=0$.

Using \eqref{kom1}, we find for any $j,k=1,...,s-1$
\begin{equation}\label{geppo}
 [E_j, \EE_k] =-\frac{1}{m_j^2}\sum_{\a,\b\in
             I_j}2ip_{\a\bb}z^\a\zz^\b \d_{jk} T
            = -2i\abs{z_j}^{2m_j} \d_{jk}T,
\end{equation}
because $\sum_{\a\in I_j}p_\a z^\a =m_j\abs{z_j}^{2m_j}$.
Thus we have
\[
 f_*T= \frac{i}{2\abs{z_j}^{2 m_j}}f_*[E_j, \EE_j]=
 \frac{i}{2\abs{z_j}^{2 m_j}}[f_*E_j, f_*\EE_j] \qquad\text{for all $j=1,\dots,s-1$.}
\]
Notice the commutation relations
\begin{equation}\label{oup}
\begin{split}
 [E_j, Z_\aa]&=0\qquad\text{if $\a \in I_s$ and   $j\le s-1$;\quad and }
\\ [Z_\a, Z_\bb]&=-2 i \d_{\a\b}T \qquad\text{if $\a,\b \in I_s$.}
\end{split}
\end{equation}
In view of \eqref{oup} and \eqref{geppo}, we claim that
  for any vector field $Z\in  \E\oplus\W_s $  there exist a real function
 $\s $ and  a vector field $U\in  \E\oplus\W_s $ such that
\begin{equation}\label{ghigo}
  [Z, \ZZ]=i\sigma T+ U-\UU.
 \end{equation}
Formula \eqref{ghigo} can be checked by a routine computation.
Here and hereafter, with slight abuse of notation we denote sections of a
bundle with the same notation of the bundle.
 The  claim applies
to $Z = f_*(E_j) \in  \E\oplus \W_s $, for $j=1,...,s-1$.
The vector fields $f_* W \in \W_1\oplus...\oplus \W_{s-1}$ and $f_* T$ are then orthogonal  with respect to the Levi form and the proof is concluded.
\end{proof}

\section{CR mappings in generalized ellipsoids}
\setcounter{equation}{0}
\label{SFRT}
\subsection{Computation of the CR factor}

\begin{lemma}
 \label{SFR}
 Let $N\subset M$ be an open set and let $v$ be a CR function in $N$ such that
$v_{,\a\b}=0$ for all
$\a,\b =1,...,n$. Then for any $\a \in I_j$, $j=1,...,s$, we have
\begin{equation}
 \label{WW}
 Z_\a T v =\frac{ n_j  (m_j-1)}{2im_j ( n+1) |z_{j}|^{2m_j}} W_\a v,
\end{equation}
where $W_\a = Q(Z_\a)$ is the vector field \eqref{sfer}.
\end{lemma}

\begin{proof}
We use the third order commutation formulae in \cite[eq.~(2.1)]{JL}. Because
$v_{,\g\a}=0$, we have
\begin{equation*}
\begin{aligned}
 v_{,\g\bb\a}= v_{,\g\a\bb}-i h_{\a\bb}v_{,\g 0} - R_\g{}^\r{}_{\a\bb}v_{,\r} =
 -i h_{\a\bb}v_{,\g 0} - R_\g{}^\r{}_{\a\bb}v_{,\r}.
\end{aligned}
\end{equation*}
On the other hand, by the commutation formula $v_{,\g\bb}- v_{,\bb\g} =
i h_{\g\bb}v_{,0}$ and since $v_{,\bb}=0$, we get $
 i h_{\g\bb}v_{,0\a}= -i h_{\a\bb}v_{,\g 0} -  R_\g{}^\r{}_{\a\bb}v_{,\r}
$. Contracting with $h^{\g\bb}$ yields
\[\begin{aligned}
 i (n+1) v_{,0\a} & = -     R_\g{}^\r{}_{\a}{}^\g v_{,\r}
= -R_\g{}^\g{}_\a{}^\r v_{,\r} =-R_\a{}^\r v_{,\r}  =
\frac{n_j(m_j-1)}{ 2m_j\abs{z_j}^{2m_j}} Q_\a^\r v_{,\r}
  \end{aligned}
\]
and the proof is concluded.
\end{proof}

\begin{theorem}
  \label{oneone}
Let  $N\subset M$ be a connected open set and let $f:N\to
f(N)\subset M$ be a CR diffeomorphism with CR
factor $\la =u^{-1}$. Then, either $u$ is a constant or there
exist $k\in\R\setminus\{ 0\}$ and  $(a_s, a^{n+1})= (a_s, t_0+i \abs{a_s}^2)\in \C^{n_s +1} $  such that
 $u=  k^{2}
\abs{z^{n+1}+ a^{n+1}+2iz_s\cdot \bar a_s}^{2}$.\end{theorem}

\begin{proof}The argument here is
similar to \cite{JL}.
The torsion $A_\theta$ of $\vartheta$ vanishes, as noted in \eqref{tower}, $A_{\a\b}=0$.
On the other hand, denoting by
$\widetilde A = A_{f^* \vartheta}$ the torsion of $\widetilde \vartheta =f^*\vartheta$, we have
$
    A_{f^*\vartheta}(Z,W) = A_\theta(f_* Z,f_*W)=0,$
 for all $ Z,W\in T^{1,0}N.
$
 Thus we also have $\widetilde A_{\a\b}=0$.
From \eqref{tors}, we deduce that $u$ satisfies the system of
equations
    $ u_{,\a\b}=0. $

By Theorem \ref{bad}, the assumptions of Proposition \ref{ricc} hold and therefore $f$ preserves
the Ricci curvature of $\vartheta$. By \eqref{cari}, we have the system of equations
\begin{equation}
 \label{sesto}
 (n+2)\Big\{
  u_{,\a\bar\b}+u_{,\bar\b\a} -\frac{2}{u} u_{,\a} u_{,\bar\b}\Big\}
  +
  \Big\{
     \Delta u - \frac{2(n+2)}{u} |\D u|^2\Big\} h_{\a\bar\b}=0.
\end{equation}
Then, the function $w =\log u$ satisfies the system of equations
$ w_{,\a\bb}+ w_{,\bb\a} = \frac{\Delta w}{n} h_{\a\bb}$. By  \cite[Proposition
3.3]{L}, $w$ is locally the real part of a CR  function $F $,
i.e., we have locally
\[
    u =e^{(F+\ol F)/2}= v\bar v = |v|^2 ,
\]
where  $ v :=e ^{{F}/{2}}$ is a CR function.
Since $u_{,\a\b}=0$, we also have $v_{,\a\b}=0$.

The function $ g:  = T v $ satisfies $g_{,\bar \a} = Z_\aa T v = T
Z_\aa v  = 0$. (Recall that $T$ and $Z_\aa$ commute.) By Proposition
\ref{radioso}, we have $W_\a u=0$ for all $\a \in I_1\cup...\cup
I_{s-1}$. This implies $W_\a v = 0$ for the same indexes. By
Lemma \ref{SFR}, we deduce that for any $\a=1,...,n$ we also have
$
    g_{,\a} =  Z_\a T v =0$. The equations $g_{,\a}=g_{,\aa}=0$ imply that
$g$ is locally constant. Therefore there exist a constant $k\in \C$ and a
function $\psi=\psi(z,\bar z)$
such that $v(z, \zz, t)= kt + \psi(z, \zz)$. Since $v$ is CR, it must be
$\p_\aa \psi - i k \p_\aa p=0$, which means
\[
      v  = k\big(t +i p(z,\bar z)\big) + \phi(z),
\]
for some holomorphic function $\phi$. Possibly multiplying $v$ by a unitary
complex number, we can assume that $k$ is real.

Moreover, for any $\a\in I_1\cup\dots\cup I_{s-1}$ we have
$ W_\a \phi  =  0$, because $W_\a v = 0$.
This fact implies that $v$ depends locally only on $\abs{z_1},\dots
,\abs{z_{s-1}} $ and, if $I_s\neq\varnothing$, on $z_s$.
From $v_{,\a\b} = 0$ and since $\D_\a Z_\b=0$ for all $\a,\b\in I_s$, see \eqref{lixor}, we deduce that  for all $\a,\b\in I_s$ we have
$
    \frac{\p^2\phi}{\p z^\a \p z^\b } = 0,
$
which finally yields
\begin{equation}
  \label{viva}
   v(z, \zz, t) = k (t+ i p(z, \zz) )+ z_s\cdot \bar d +c
\end{equation}
for some $d\in \C^{n_s}$ and $c\in\C$.
In order to find the imaginary part of $c$, we observe that if  $u$ solves
the system \eqref{sesto}, it also solves the contracted equation $\Delta u -
\frac{n+2}{u}\abs{\nabla u}^2 =0$. After a computation, this implies that  $v$
solves the
equation
\begin{equation}
 \label{vuva}
   i \big( \bar v v_{,0} - v \bar v_{,0}\big) = h^{\ga\bar\mu} v_{,\ga}
   \overline{v_{,\mu}}.
\end{equation}
Plugging  \eqref{viva} into \eqref{vuva}, we get
$
  k \mathrm{Im}( c) = {|d|^2}/ {4}.
$
If $k=0$ then $d=0$, $\phi(z)=c$  and $v$ is
 constant.
If  $k\neq 0$ then, letting
 $t_0 =   \Re (c)/k$,  we obtain
$  v = k \Big \{t+ t_0 +
    z_s\cdot  \bar d/k + i  \big( p +|d|^2/(4k^2) \big )\Big\} .
$ Letting
 $\bar d/k=2i \bar a_s$, we get
\begin{equation}
  \label{viva2}
   v = k \Big \{t+t_0  +i p +2i z_s\cdot\bar a_s + i \abs{a_s}^2 \Big\} = k
\big\{z^{n+1} + a^{n+1} +2iz_s\cdot\bar a_s\big\},
\end{equation}
where $z^{n+1}= t+i p(z,\zz)$ and $a^{n+1} = t_0 + i \abs{a_s}^2$.
This concludes the proof.
\end{proof}

 \subsection{Levi-isometric mappings}
\label{leviso}
Let $M$ be the surface \eqref{dominio}
endowed with the pseudohermitian structure $\vartheta$ introduced in
\eqref{conto}.
We say that a CR diffeomorphism $\psi:N\to N' $ is \emph{Levi-isometric with respect to
$\vartheta$} if  $\psi^*\vartheta =
\vartheta$.
 A Levi isometric mapping $\psi$ satisfies $L(\psi_* Z, \psi_*\bar W) = L
(Z,\bar W)$ for all $Z,W\in T^{1,0}N$.
Moreover, since we trivially have
$
 \psi_* (T^{\psi^*\theta}) =(T^\theta)_\psi
$,   it turns out that a Levi isometric mapping satisfies
\begin{equation}
\label{acci}
\psi_*T=T_\psi. \end{equation}

\begin{theorem}  \label{isobar}
Let $N,N'\subset M$ be connected open sets and let $\psi:N\to N' $
be a Levi-isometric mapping with respect to $\vartheta$. Then, there
exists  a permutation $\sigma$ of  $\{1, \dots, s-1\}$ such that
              $m_{\s(j)}= m_{j}$ and  $n_{\s(j)}= n_{j}$ for any $j=1, \dots,
              s-1$, there are
 unitary matrices  $B_j\in U(n_j)$, and,   if $n_s\ge 1$, there are    $B_s\in
U(n_s)$ and
          a   vector $(b_s, b^{n+1})= (b_s, t_0+ i\abs{b_s}^2)\in \C^{n_s}\times\C$  such that for all
          $(z,t)\in N$ we have
\begin{equation}
  \label{LeviU}
 \psi(z,z^{n+1})=\left( B_1 z_{\sigma(1)}, \dots, B_{s-1}z_{\sigma(s-1)}, B_sz_s
+ b_s,
  b^{n+1}+ z^{n+1}+2i  ( B_sz_s \cdot  \bar b_s)\right).
\end{equation}
\end{theorem}

 We start with an easy lemma.

\begin{lemma}
 \label{ORTO}
Let $D\subset\C^d$, $d\geq 2$, be an open connected set and let $\zeta
:D\to\C^d$ be a nonconstant   holomorphic mapping such that $|\z(z)|$ is
constant if $|z|$ is constant, for $z \in D$. Then there exists $B\in GL(d,\C)$
such that $\z(z) = B z$ and $B^* B = \rho^2 I$ for some  $\rho >0$.
\end{lemma}

\begin{proof}
Assume without loss of generality that there exists $z\in D$ such that $|\z(z)|=|z|=1$. This can be achieved multiplying $\zeta$ by a positive constant.
Then, by the Poincar\'e--Alexander theorem, see \cite{Alexander,Tanaka62,Rudin81},  $\z$ is the restriction of an automorphism of the unit ball $ B_1 :  =\{z\in\C^d:|z|<1\}$.
  Thus, see \cite{Krantz,Rudin81},  there exist a unitary matrix $B\in U(d)$ and
$a\in\C^d$ with $|a|<1$ such that $\z (z) = B\phi_a(z)$ for all
$z\in D$, where $
   \phi_a(z)=\frac{a-Pz  {-}\sqrt{1-\abs{a}^2}Qz}{1- z\cdot\bar a },
$
with $Pz= {\frac{z\cdot \bar a}{\abs{a}^2}a}$ and  $Qz=z-Pz$.
When $a=0$ we have $\phi_a(z) = -z$.
If $a\neq 0$,  $\phi_a$ takes
    $b B_1$ to $b B_1$ but it does not take any  other (open piece of) sphere
      $b  B_r $ with $r\neq 1$ to a sphere centered at the origin.
Then we have $a=0$ and the Lemma follows.
\end{proof}

\begin{proof}[Proof of Theorem \ref{isobar}] All our claims along the proof are of a local nature. In  the coordinates $(z,t)\in\C^n\times\R$  on $M$,  we have  $\psi = (\z^1,...,\zeta^n,\tau)$ with $\zeta^\b: N\to \C$,  $\b=1,...,n$, and $\tau:N\to\R$.
We first notice that  $Z_\aa \zeta^\b =0$ for all $\a,\b=1,...,n$,
because $\psi$ is a CR mapping. Moreover, we have $T\tau
=1$ and $T\zeta^\b= 0$ because $\psi_*T=T$, by \eqref{acci}. Then, for $j=1,
\dots,s$
the functions $\zeta_j =\z_j(z)$ are holomorphic and $\tau =
t+v(z,\bar z)$ for some real function $v$.

By Theorem \ref{bad},  there exists a permutation  $\sigma$ of $\{1,
\dots, s-1\}$ such that $\psi_* \W_j = \W_{\sigma (j)}$. In the
following we let $j'=\sigma(j)$. In particular, we have
$n_{j'}=n_{j}$ for all $j=1,...,s-1$.

Fix $j\in\{1,...,s-1\}$.
Let $V\in \W_j$ with $|V|=1$. Since $\psi$ is Levi isometric, $\psi$ preserves the sectional curvature of $\vartheta$, $k(\psi_* V) = k(V)$. By \eqref{sec}, we deduce that
\begin{equation}
 \label{SIGMA1}
   \frac{m_j-1}{m_j} |z_j|^{-2m_j} =
   \frac{m_{j'}-1}{m_{j'}} |\z_{j'}(z)|^{-2m_{j'}}.
\end{equation}
With the notation $z_j^* = (z_1,...,z_{j-1}, z_{j+1},...,z_s)$, consider for
fixed $z_j^*$ the mapping $z_j\mapsto \zeta_{j'}(z_j^*; z_j)$. This mapping is
nonconstant and holomorphic from an open subset of $\C^{n_j}$ to $\C^{n_j}$.
Moreover, by \eqref{SIGMA1} it  takes (pieces of) spheres of $\C^{n_j}$
centered at the origin into spheres centered at the origin. By Lemma \ref{ORTO},
there exist  $B_j\in GL(n_j,\C)$  and $\rho_j>0$ such that $\zeta_{j'}(z) = B_j (z_j^*)z_{j}$ with  $B_j^*B_j =\rho_j^ 2  I$. Here $B_j=B(z_j^*)$
 is holomorphic, while   $\rho_j=\rho_j(z_j^*, \bar z_j^*)$.
Therefore,
\eqref{SIGMA1} becomes
\begin{equation*}
 \abs{z_j}^{2(m_j-m_{j'})}
 = \frac{m_j-1}{m_j}\, \frac{m_{j'}}{m_{j'}-1} \rho_j(z_j^*, \bar z_j^*)^{2m_{j'}}.
\end{equation*}
Both the left-hand side and the right-hand side must be constant. Therefore $m_j=m_{j'}$ and
$\rho (z_j^*, \bar z_j^*)=1$. Ultimately we have for any $j\le s-1$,  $\z_{j'}(z)= B_j z_j$, for some constant matrix  $B_j\in U(n_j)$.

Next we claim that, if $n_s\ge 1$ then $\z_s$ depends only on $z_s$.
To prove the claim it suffices to show that for all $j\le s-1$ we have
\begin{subequations}
 \begin{align}
 \label{duedue}W_\l \z^\g & =0 \qquad\text{for all $\l\in I_j$, $\g\in
I_s$;\quad and}
\\ \label{tretre}
E_j\z^\g &  =0 \qquad \text{for all   $\gamma\in I_s.$}
\end{align}
\end{subequations}
To prove \eqref{duedue}, fix $j\le s-1$ and $\l\in I_j$. Then
\begin{equation*}
\begin{aligned} \psi_*W_\l & = \sum_{k=1}^{s-1} \sum_{\g\in I_{k}} (W_\l \z^\g)(\p_\g)_\psi +
\sum_{\g\in I_s}  (W_\l \z^\g)(\p_\g)_\psi +(W_\l \t) (\p_t)_\psi  \qquad \text{(since $\psi$ is CR)}
\\& =  \sum_{k=1}^{s-1} \sum_{\g\in I_{k}} (W_\l \z^\g)(Z_\g)_\psi +
\sum_{\g\in I_s}  (W_\l \z^\g)(Z_\g)_\psi .
\end{aligned}
\end{equation*}
But $\psi_*W_\l\in \W_{j'} $. Then all the terms in the last sum  vanish and \eqref{duedue} is proved.

To prove \eqref{tretre} start by computing $\psi_*E_j$ for $j\le s-1$:
 \[\begin{aligned}
         \psi_*E_j&= \sum_{k'=1}^{s-1}
 \sum_{\substack{\g\in I_{k'}\\ \mu\in I_k}} B_{\g\m} E_j z^\m (Z_\g)_\psi
         +\sum_{\g\in I_s} E_j\z^\g (Z_\g)_\psi
\\&= \sum_{\substack{\g\in I_{j'}\\ \mu\in I_j}} B_{\g\m} \frac{1}{m_j	}z^\m (Z_\g)_\psi +
\sum_{\g\in I_s} E_j\z^\g (Z_\g)_\psi
= (E_{j'})_\psi + \sum_{\g\in I_s} E_j\z^\g (Z_\g)_\psi.
         \end{aligned}
 \]
 Taking the Levi-length,
 we find $
\abs{\psi_* E_j}^2 =\abs{(E_{j'})_\psi}^2 + \sum_{\g,\r\in I_s}(E_j\z^\g) (\overline E_j \bar\z^\r)( h_{\g\bar\rho})_\psi.
 $
But we have $\abs{\psi_* E_j}^2=\abs{E_j}^2 = 2\abs{z_j}^{2m_j}$ and $\abs{(E_{j'})_\psi}^2=2\abs{\z_{j'}}^{2m_j}=
2 \abs{z_j}^{2m_j}$. Moreover it is $( h_{\g\bar\rho})_\psi=2\d_{\g\r}$, because $\g,\rho\in I_s$.
 Then we conclude  that $E_j\z^\g=0$, as required.

 Next, we compute $v$ and $\z_s$.
  Let $\a\in I_j$, where $j\le s-1$. Recall that $\tau= t+ v(z, \zz)$. Since
$\p_\a\z_s=0$, we have
  \[
\begin{aligned}
 \psi_*Z_\a & =\sum_{\b\in I_{j'}}B_{\b\a}(\p_\b)_\psi + \big(\p_\a v +i m_j\abs{z_j}^{2(m_j-1)}\zz^\a \big)(\p_t)_\psi
\qquad \text{
(since $\psi$ is CR) }
\\
 & = \sum_{\b\in I_{j'}}B_{\b\a}(Z_\b)_\psi =
\sum_{\b\in I_{j'}}B_{\b\a} \big\{(\p_\b)_\psi + i m_{j'}\abs{\z_{j'}}^{2(m_{j'}-1)}\bar\zeta^\b(\p_t)_\psi\big\}
\\&=\sum_{\b\in I_{j'}}B_{\b\a}(\p_\b)_\psi + i m_{j}\abs{z_{j}}^{2(m_{j}-1)}\bar{z}^\a(\p_t)_\psi,
\end{aligned}
\]
where we used $\abs{\z_{j'}}=\abs{z_j}$, $m_{j'}=m_j$ and $\sum_{\b\in
I_{j'}}B_{\b\a}\overline{B_{\b\g}}= \d_{\a\g}$, if $\a, \g\in I_j$. Comparing
the first and third lines we get $v_\a=0$. Thus $v$ depends only on
$z_s, t$.

Finally, we find $\zeta_s$ and $v$  when $n_s\geq 1$.
For  $\a\in I_s$ we have
\begin{equation}
\begin{aligned}\label{8000}
 \psi_*Z_\a &= \sum_{\b\in I_s}(\p_\a \z^\b) (\p_\b)_\psi
  + \big\{i\bar z^\a +\p_\a v \big\} T_\psi,
\end{aligned}
\end{equation}
as well as
\begin{equation}\label{ciop2}
\begin{aligned}
 \psi_*Z_\a &= \sum_{\b\in I_s}  (\p_\a\z^\b )(Z_\b)_\psi
 =  \sum_{\beta\in I_s} (\partial_\a \zeta^\b )\big(( \partial _\b)_\psi + i \bar \zeta^\b T_\psi\big).
\end{aligned}
\end{equation}
Since $\psi$ is Levi isometric, we have  $2\d_{\a\g}=  L(\psi_*Z_\a,
\psi_*Z_\gg) $
for $\a, \g\in I_s$. Using   formula \eqref{ciop2}, we obtain
$ \sum_{\b\in I_s}(\p_\a \z^\b)(\p_\gg\bar  \zeta^\b)=\d_{\a\g}
$, for all $\a, \g\in I_s$.
Therefore it must be $ \z_s(z) = b_s+Bz_s$, for some   $B\in U(n_s)$ and $b_s\in\C^{n_s}$.
Moreover, comparing the coefficients of $T$ in \eqref{8000} and \eqref{ciop2},
 we obtain the equation
$
   i \bar z^\a +\partial_\a  v = i \sum_{\b\in I_s} \bar\zeta^\b \partial _\a \zeta^\b,
$
which implies  $  \p_\a v= {i}\sum_{\b\in I_s} B_{\b\a} \bar b^\b
$.
Since $v$ is a real function, we finally find $ v = t_0 -
2\mathrm{Im}(Bz_s\cdot \bar b_s)$ for some $t_0\in\R$.

The structure \eqref{LeviU} of the isometry $\psi$ is now determined locally.
The proof is concluded because $N$ is connected.
\end{proof}

\section{Proof of Theorem \ref{bad}}
\setcounter{equation}{0}
\label{brutta}
This section is devoted to the proof  of  Theorem \ref{bad}. The proof is rather
involved, but we were not able to find a more direct one.  In many situations,
the study of this case  can be
 avoided
 for trivial  dimensional reasons, see the discussion before the statement of
Theorem \ref{bad}.

 \begin{proposition}
 \label{biba}
  Let $N\subset M$ be an open set and let $f:N\to
f(N)\subset M$ be a CR diffeomorphism such that
$f_*(\E\oplus\W_s)=\mathcal W_j$ for some $j=1,...,s-1$. Then
there exists $\nu\in\N\cup \{0\}$ such that
\begin{equation}\label{cairo}
         M= \big \{ (z_1, z_2, z^{n+1})\in\C^{\nu+2}\times\C^\nu\times\C:\Im
z^{n+1}
=|z_1 |^{2m_1} + |z_2|^2, \,\text{ and $z_1\neq 0$} \big\}.
\end{equation}
Moreover, we have  $\l R\circ f= R$ where $\l$ is the CR factor of $f$.
\end{proposition}

\begin{proof}
For some $k=1,...,s-1$ we have $f_*(\mathcal W_k)=\E\oplus\W_s$. Since $f$ is a CR diffeomorphism, it must be
$\dim \W_j = \dim(\E\oplus\W_s)=\dim \W_k$. In other words,
\begin{equation}
n_j= n_s+s = n_k \label{picopaco}
\end{equation}
 For any
$V\in\mathcal W_k$ with $|V|=1$ we evaluate  $S(V): = S(V,\VV,V,\VV)$. By \eqref{web},
 \eqref{detto},   \eqref{ture},
and \eqref{Q2} we get
\begin{equation}
 \label{wetzen}
\begin{split}
   S(V) &
  =\frac{1}{n+2}\left\{ (2n_k-n-2)\frac{m_k-1}{m_k }   |z_k|^{-2 m_k}
   +\frac{2 R}{n+1}\right\}.
\end{split}
\end{equation}
Since  $f_*V\in\E\oplus\W_s$, all the terms involving curvature in
$S(f_* V)$ vanish and we get
\begin{equation}
 \label{messer}
\begin{split}
   S(f_* V)
   & = \frac{2 R\circ f }{(n+1)(n+2)}|f_* V|^4
    = \la^2 \frac{2 R\circ f }{(n+1)(n+2)}.
\end{split}
\end{equation}
By the CR invariance  \eqref{chevi}, we deduce from
\eqref{wetzen} and \eqref{messer}
\begin{equation}
 \label{secolano}
  \frac{2\la R\circ f }{n+1} = (2n_k-n-2)\frac{m_k-1}{m_k }   |z_k|^{-2 m_k}
   +\frac{2 R}{n+1}.
\end{equation}

 Let $Z\in\E\oplus\W_s$ with $|Z|=1$.
 We have
\begin{equation}
 \label{loeffel}
   S(Z) = \frac{2R}{(n+1)(n+2)}.
\end{equation}
Since $f_* Z\in\mathcal W_j$ for some $j\le s-1$,
arguing as in \eqref{wetzen} we find
\begin{equation}
 \label{fork}
 S(f_*Z) = \frac{\la^2}{n+2}\left\{ (2n_j-n-2)\frac{m_j-1}{m_j }
 |\z_j|^{-2 m_j}
   +\frac{2 R\circ f}{n+1}\right\},
\end{equation}
where  we let $(\zeta,\z^{n+1}) = f(z,z^{n+1})\in M$. By the CR invariance
\eqref{chevi}, we obtain
\begin{equation}
 \label{pup}
   \la \left\{ (2n_j-n-2)\frac{m_j-1}{m_j }
 |\z_j|^{-2 m_j}
   +\frac{2 R\circ f}{n+1}\right\} = \frac{2R}{n+1}.
\end{equation}
Comparing \eqref{secolano} and \eqref{pup} we get
\begin{equation}
 \label{putto}
 (2n_k-n-2)\frac{m_k-1}{m_k }   |z_k|^{-2 m_k}
  +  \la (2n_j-n-2)\frac{m_j-1}{m_j }
 |\z_j|^{-2 m_j}
   =0 .
\end{equation}
Recall that by \eqref{picopaco}, it must be
$n_j=n_k$.
Therefore, \eqref{putto} becomes
\begin{equation*}
 (2n_j-n-2)\Big\{ \frac{m_k-1}{m_k }   |z_k|^{-2 m_k}
  +  \la  \frac{m_j-1}{m_j }
 |\z_j|^{-2 m_j}\Big\}
   =0 .
\end{equation*}
But the curly bracket is positive. Then, we have  $n_j=n_k=(n+2)/2$. Moreover,
it
must be $j=k$, because if $j\neq k$ the condition
$n_j+n_k\leq n$ is not satisfied.
Finally, using
$n_i\ge 2$ for all $i\le s-1$ and  $n_s = n_j-s$,
we get
\[
 n= n_1+...+n_s=n_j+n_s+\sum_{i\neq j} n_i
\ge n_j+n_s +2(s-2)= 2n_j+s-4=
n+s-2.
\]
This gives $s\le 2$. If $s=1$, then $M$ is the surface   $\Im(z^3)=
(\abs{z^1}^2+ \abs{z^2}^2)^{m_1}$. If    $s=2$, we have $n=n_1+n_2= n/2 +
1+n_2$,
which implies $n_1=n_2+2$ and the domain has the form  \eqref{cairo}.
\end{proof}

\begin{proposition}
  \label{uai}
Let $M$ be as in \eqref{cairo} and   let $f:N\to N' $   be a
CR diffeomorpshism
  such that    $f_*(\E\oplus\W_2 ) = \mathcal
  W_1$.  Then the CR  factor $\la$ of $f$ satisfies
   $E_1\la = -\la$.
\end{proposition}

\begin{proof}
Fix $\s,\mu\in I_1$ and  let $W= \zz^\s Z_\m-\zz^\m Z_\s=
\zz^\s\p_\m-\zz^\m\p_\s$. Notice that $L(W,\bar E_1)=0$ and thus
$W\in\W_1$. We also have $
 [W, \WW]= z^\s\p_\s +z^\m\p_\m - \zz^\s\p_\ss -\zz^\m\p_\mm$ and
\begin{equation}
 \label{pralleno}
 [[W, \WW], \EE_1]=\Big[z^\s \p_\s + z^\m \p_\m - \zz^\s\p_\ss - \zz^\m \p_\mm, \frac{1}{m_1} \sum_{\b\in I_1}\zz^\b\p_\bb - i\abs{z_1}^{2m_1}\p_t\Big]=0.
\end{equation}
Finally, we have
$
i \vartheta([W,\bar W]) = -id\vartheta(W,\bar W) =
 2m_1\abs{z_1}^{2(m_1-1)}\big(\abs{z^\s}^2+\abs{z^\m}^2\big).
 $

Since   $f_*W \in \E\oplus \W_2$, as in the proof of Proposition
\ref{radioso}, see \eqref{ghigo}, we have
  \begin{equation}
  \label{este2}
  f_*[W, \WW]= [f_* W, f_*\WW]= F -\bar F + i k T ,
  \end{equation}
for some $F\in   \mathcal E {\oplus \W_2}$ and some real function $k$ on
$f(N)$. Since
$f_*E_1\in   \W_1$,
$f_*[W,\WW]$ and $f_* \bar E_1$
are orthogonal by \eqref{este2}. Then, also using \eqref{pralleno}, we get
\begin{equation}
\label{sposto}
 \begin{split}
 0= d\vartheta(f_*[W, \WW], f_*\EE_1)
    =- (f_*\EE_1)(\vartheta(f_*[W, \WW])),
 \end{split}
\end{equation}
that is equivalent with
$
 E_1 \big( \lambda \abs{z_1}^{ {2(m_1-1)}}(\abs{z^\s}^2 +\abs{z^\m}^2)\big)=0.
$
Since $\s,\m$ are arbitrary, this
implies $E_1 \big( \lambda \abs{z_1}^{2m_1}\big)=0$ and eventually
$E_1\l+\l=0$, because $E_1  \abs{z_1}^{2m_1} =  \abs{z_1}^{2m_1}$.
\end{proof}

\bigskip

\begin{proof}[Proof of Theorem \ref{bad}]
Assume by contradiction that $f_*(\E\oplus\W_s) = \W_j$ for some
$j=1,...,s-1$. Then, by Proposition \ref{biba}, $M$ is of the form
\eqref{cairo}, and by Proposition \ref{uai} we have $E_1\l=-\l$,
where  $\la =u^{-1}$ is the CR factor of $f$. We have $f_*
\W_1 = \W_2\oplus \E$ and $f_* (\W_2\oplus \E) = \W_1$.

In terms of $u$ we have $E_1 u = u$. Note that \eqref{D1} implies that
  $\D_{\EE_1} E_1=0$. Thus
\begin{equation*}
   \D^2u(E_1, \EE_1)+\D^2u(\EE_1, E_1)-\frac 2u |E_1u|^2=0.
\end{equation*}
Since $\Ric(E_1, \bar E_1)=0$,    \eqref{ochetta} becomes
\begin{equation}\label{flauto}
 \begin{aligned}
 \Ric(f_*E_1, f_*\bar E_1)&=  \frac{1}{2u} \Big\{\Delta u -\frac {2(n+2)}{u}|\D u|^2
 \Big\}|E_1|^2
\\&= \frac 12 \Big\{\Delta u -\frac {2(n+2)}{u}|\D u|^2
 \Big\}|f_*E_1|^2.
 \end{aligned}
\end{equation}
On the other hand, since $f_*E_1\in \W_1$, comparing   \eqref{ture} and
\eqref{coccobill}, we get
\begin{equation}
  \label{piffero}
  \Ric(f_*E_1, f_*\bar E_1) = \frac {1}{n_1-1}( R\circ f) \left| f_* E_1\right|^2,
\end{equation}
and therefore  \eqref{flauto} becomes
\begin{equation}
 \label{ottavino}
  2\frac{R\circ f}{n_1-1} =
  \Delta u -\frac {2(n+2)}{u}|\D u|^2.
\end{equation}
By Proposition \ref{biba} we have
 $R\circ f = u R$, and from \eqref{change4}  we obtain
\begin{equation*}
u R = R\circ f=\wt R= uR+(n+1)\Big\{ \Delta u-\frac {n+2}{u}|\D u|^2\Big\},
\end{equation*}
that gives
$\displaystyle
  \Delta u=\frac {n+2}{u}|\D u|^2$.
Inserting this identity into  \eqref{ottavino} and using
formula
\eqref{coccobill}, we obtain
\begin{equation}
 \label{pink}
 \frac{|\D u|^2}{u^2}=\frac{m_1-1}{m_1(n+2) |z_1|^{2m_1}}n_1= \frac{m_1-1}{2m_1\abs{z_1}^{2m_1}},
\end{equation}
because, $n_1=\nu+2$ and  $n_2=\nu$, so that $\frac{n_1}{n+2}=\frac 12.$
On the other hand,
by
$E_1u=u$ and $\abs{E_1}^2=2\abs{z_1}^{2m_1}$
we have
$   |\D u|^2
  \ge
  \frac{|E_1 u|^2}{|E_1|^2}=\frac{u^2}{2\abs{z_1}^{2m_1}},
$ which   contradicts   \eqref{pink}. The proof is concluded.
\end{proof}



\providecommand{\bysame}{\leavevmode\hbox to3em{\hrulefill}\thinspace}
\providecommand{\MR}{\relax\ifhmode\unskip\space\fi MR }
\providecommand{\MRhref}[2]{%
  \href{http://www.ams.org/mathscinet-getitem?mr=#1}{#2}
}
\providecommand{\href}[2]{#2}

\bigskip \noindent\sc \small  Roberto Monti
\\ Dipartimento di Matematica Pura ed Applicata,
Universit\`{a} di Padova  (Italy)
\\Email: \tt   monti@math.unipd.it

\bigskip
\noindent\sc \small  Daniele Morbidelli \\ Dipartimento di Matematica, Universit\`a di Bologna (Italy)
\\Email: \tt morbidel@dm.unibo.it. \rm

\end{document}